\documentclass[11pt,reqno]{amsart}
\usepackage{amsmath,amsxtra,latexsym,amsthm,amssymb,amscd,pb-diagram}
\usepackage{mathrsfs,mathabx,amsfonts}
\usepackage[colorlinks=true,linkcolor=magenta,citecolor=blue]{hyperref}
\usepackage[margin=1in]{geometry}
\usepackage{bm}
\usepackage{color}
\usepackage{makecell,multirow,diagbox}
\usepackage{mathtools}
\usepackage[displaymath, mathlines]{lineno}

\usepackage{graphicx}
\usepackage{epsfig}
\usepackage{array}
\usepackage{enumitem}

\newtheorem{theorem}{Theorem}[section]
\newtheorem{proposition}{Proposition}[section]
\newtheorem{lemma}{Lemma}[section]

\newtheorem{remark}{Remark}[section]

\newcommand{\R}{\mathbb{R}}

\newcommand{\Z}{\mathbb{Z}}
\newcommand{\N}{\mathbb{N}}
\newcommand{\e}{\varepsilon}
\newcommand{\ity}{\infty}
\newcommand{\dps}{\displaystyle}
\newcommand{\f}{\displaystyle\frac}

\begin{document}
\title[A system of semi-linear $\sigma$-evolution equations with double dissipation]{On the Cauchy problem for a weakly coupled system of semi-linear $\sigma$-evolution equations with double dissipation}
\subjclass{35A01, 35B44, 35L56}
\keywords{Weakly coupled system, Double damping, Global solutions, Blow-up, Lifespan}
\thanks{\textit{Corresponding author:} Tuan Anh Dao}

\maketitle
\centerline{\scshape Yingli Qiao}
\medskip
{\footnotesize
	\centerline{School of Mathematical Sciences, Nanjing Normal University}
	\centerline{Nanjing, 210023, China}
	\centerline{Email: 210902076@njnu.edu.cn}}
\medskip

\centerline{\scshape Tuan Anh Dao}
\medskip
{\footnotesize
	\centerline{Faculty of Mathematics and Informatics, Hanoi University of Science and Technology}
	\centerline{No.1 Dai Co Viet road, Hanoi, Vietnam}
	\centerline{Email: anh.daotuan@hust.edu.vn}}

\begin{abstract}
In this paper, we would like to consider the Cauchy problem for a multi-component weakly coupled system of semi-linear $\sigma$-evolution equations with double dissipation for any $\sigma\ge 1$. The first main purpose is to obtain the global (in time) existence of small data solutions in the supercritical condition by assuming additional $L^1$ regularity for the initial data and using multi-loss of decay wisely. For the second main one, we are interested in establishing the blow-up results together with sharp estimates for lifespan of solutions in the subcritical case. The proof is based on a contradiction argument with the help of modified test functions to derive the upper bound estimates. Finally, we succeed in catching the lower bound estimate by constructing Sobolev spaces with the time-dependent weighted functions of polynomial type in their corresponding norms.
\end{abstract}

% \linenumbers
\tableofcontents

%==========================================================
\section{Introduction}
Let us consider the following Cauchy problem for weakly coupled system of semi-linear $\sigma$-evolution equations with two dissipative terms:
\begin{equation} \label{MainSytem}
\begin{cases}
\partial^2_t u_1+ (-\Delta)^\sigma u_1+ \partial_t u_1+ (-\Delta)^{\sigma} \partial_t u_1= |u_{\text{k}}|^{p_1}, &\quad x\in \R^n,\, t \ge 0 \\
\partial^2_t u_2+ (-\Delta)^\sigma u_2+ \partial_t u_2+ (-\Delta)^{\sigma} \partial_t u_2= |u_1|^{p_2}, &\quad x\in \R^n,\, t \ge 0 \\
\quad \vdots \\
\partial^2_t u_{\mathtt{k}}+ (-\Delta)^\sigma u_{\mathtt{k}}+ \partial_t u_{\mathtt{k}}+ (-\Delta)^{\sigma} \partial_t u_{\mathtt{k}}= |u_{\mathtt{k}-1}|^{p_{\mathtt{k}}}, &\quad x\in \R^n,\, t \ge 0 \\
u_\ell(0,x)= u_{0\ell}(x),\quad \partial_t u_\ell(0,x)= u_{1\ell}(x), &\quad x\in \R^n,\,\ell=1,2,\cdots,\mathtt{k},
\end{cases}
\end{equation}
where $\sigma\ge 1$ is assumed to be any fractional number and $\mathtt{k}\ge 2$. With $\ell=1,2,\cdots,\mathtt{k}$, we denote by $p_\ell>1$, the power exponents of nonlinear terms, moreover, the dissipative terms $\partial_t u_\ell$ and $(-\Delta)^{\sigma} \partial_t u_\ell$ stand for the frictional dissipation (or external damping) and visco-elastic dissipation (or strong damping), respectively. The corresponding linear equation of \eqref{MainSytem} we have in mind is
\begin{equation} \label{LinearEq}
\begin{cases}
u_{tt}+ (-\Delta)^\sigma u+ u_t+ (-\Delta)^{\sigma} u_t= 0, &\quad x\in \R^n,\, t \ge 0 \\
u(0,x)= u_0(x),\quad u_t(0,x)=u_1(x), &\quad x\in \R^n.
\end{cases}
\end{equation}

First of all, let us begin with the most typical important problems for \eqref{LinearEq} and \eqref{MainSytem} with $k=1$ and $\sigma=1$, the so-called wave equations with frictional damping and visco-elastic damping, namely
\begin{equation} \label{LinearEq_sigma1}
\begin{cases}
u_{tt}- \Delta u+ u_t- \Delta u_t= 0, &\quad x\in \R^n,\, t \ge 0 \\
u(0,x)= u_0(x),\quad u_t(0,x)=u_1(x), &\quad x\in \R^n
\end{cases}
\end{equation}
and
\begin{equation} \label{MainEq_sigma1}
\begin{cases}
u_{tt}- \Delta u+ u_t- \Delta u_t= |u|^p, &\quad x\in \R^n,\, t \ge 0 \\
u(0,x)= u_0(x),\quad u_t(0,x)=u_1(x), &\quad x\in \R^n.
\end{cases}
\end{equation}
Concerning the linear equation \eqref{LinearEq_sigma1} the authors in \cite{IkehataSawada2016} succeeded in obtaining the asymptotic profile of solutions in the $L^2$ norm by assuming the weighted $L^{1,1}$ initial data. By the means of studying of asymptotic profile as $t\to \ity$, they indicated that the influence of the frictional damping is really more dominant than that of the visco-elastic one in a comparison between such two damping types. Roughly speaking, the asymptotic profile of solutions to \eqref{LinearEq_sigma1} is same as that to the following heat equation:
\begin{equation} \label{HeatEq}
\begin{cases}
v_t- \Delta v= 0, &\quad x\in \R^n,\, t \ge 0 \\
v(0,x)= v_0(x), &\quad x\in \R^n
\end{cases}
\end{equation}
for a suitable choice of data $v_0$. Furthermore, under more heavy moment conditions on the initial data for any space dimensions, some higher-order asymptotic expansions of solutions to \eqref{LinearEq_sigma1} were established in \cite{IkehataMichihisa2019} by employing Taylor expansion method effectively (see more \cite{Michihisa2021,Takeda2015}). For the treatment of the semi-linear equation \eqref{MainEq_sigma1}, the author in \cite{DAbbicco2017} proved the global (in time) solutions for any space dimensions by deriving energy estimates together with $L^1-L^1$ estimates. Telling about the influence of the two damping types as mentioned in \cite{IkehataSawada2016} to \eqref{LinearEq_sigma1}, the authors in \cite{IkehataTakeda2017} claimed that this effect is still valid for \eqref{MainEq_sigma1} with respect to the \textit{critical exponent}. Here the critical exponent is understood as a threshold to classify between global (in time) existence of small data solutions (stability of the zero solution) and finite time blow-up result of (even) small data solutions. More precisely, the critical exponent for \eqref{MainEq_sigma1} coincides with that for the semi-linear heat equation of \eqref{HeatEq} with nonlinearity term $|v|^p$ and the semi-linear classical damped wave equations with nonlinearity term $|u|^p$ as well. As usual, we call it the well-known Fujita exponent $p_{\rm Fuj}=1+2/n$. Among other things in \cite{IkehataTakeda2017}, both the existence of global (in time) solutions to \eqref{MainEq_sigma1} and the large time behavior of these global solutions have been discussed in low dimensional cases. Quite recently, the authors in \cite{MezadekMezadekReissig2020} not only have demonstrated the global (in time) existence of Sobolev solutions to \eqref{MainEq_sigma1} with suitable regularity but also investigated the interplay between the regularity for the initial data and the admissile range of exponents $p$ on these global results.

Secondly, let us continue the Cauchy problem for semi-linear $\sigma$-evolution equations with frictional and visco-elastic damping term, i.e. \eqref{MainSytem} with $k=1$ and $\sigma>1$ as follows:
 \begin{equation}\label{MainEq_sigma}
 \begin{cases}
 u_{tt}+(-\Delta)^\sigma u+u_t+(-\Delta)^\sigma u_t=|u|^p, &\quad x\in \R^n,\, t \ge 0 \\
u(0,x)= u_0(x),\quad u_t(0,x)=u_1(x), &\quad x\in \R^n.
 \end{cases}
 \end{equation}
  The relation between the two types of damping terms has been explored in \cite{DaoMichihisa2020}. The authors explained that the frictional damping is still dominant compared with visco-elastic damping for any $\sigma> 1$. The point worth noticing is that they indicated the diffusion phenomenon of $L^p-L^q$ theory with $1\le p\le q\le \infty$ under considering the corresponding anomalous diffusion problem of (\ref{HeatEq}) in the form
\begin{equation} \label{HeatEq_sigma}
\begin{cases}
v_t + (-\Delta)^\sigma v= 0, &\quad x\in \R^n,\, t \ge 0 \\
v(0,x)= v_0(x), &\quad x\in \R^n
\end{cases}
\end{equation}
with a suitable choice of data $v_0$. Moreover, they found that the frictional damping affects the asymptotic profile of the solution to (\ref{MainEq_sigma}) at small frequencies, whereas the asymptotic profile is modified by visco-elastic damping at large frequencies. The higher-order asymptotic expansion of the difference between solutions to (\ref{LinearEq}) and (\ref{HeatEq_sigma}) was also demonstrated there. By assuming additional $L^m$ regularity for the initial value, the authors in \cite{MezadekMezadekReissig2022} obtained $(L^m\cap L^q)-L^q$ estimates for Sobolev solutions on the basis of $L^q-L^q$ estimates to (\ref{MainEq_sigma}), where $q\in (1,\ity)$ and $m\in[1,q)$. After it is proved that the Sobolev solutions exist globally under selecting a suitable regularity of the initial data, the relation between the regularity hypothesis of the initial data and the nonlinear power exponent $p$ has been also reflected in \cite{MezadekMezadekReissig2022}. When $\sigma$ is an integer, the critical exponent of (\ref{MainEq_sigma}) is given by $p_{\rm crit}=1+2\sigma/n$, which is guaranteed in \cite{DaoMichihisa2020} by a result of blow-up solutions. The quantity is sometimes called the so-called modified Fujita exponent.

  During the last decades, many mathematicians have devoted their works to the study of the Cauchy problem for weakly coupled systems instead of only attacking single equations because of their wide application in various disciplines such as physics, chemistry, mechanics and so on. For this reason, let us sketch out several recent results involving the following weakly coupled systems:
\begin{equation}\label{System_sigma}
\begin{cases}
u_{tt}+(-\Delta)^\sigma u+u_t+(-\Delta)^\sigma u_t=|v|^p, &\quad x\in \R^n,\, t \ge 0 \\
v_{tt}+(-\Delta)^\sigma v+v_t+(-\Delta)^\sigma v_t=|u|^p, &\quad x\in \R^n,\, t \ge 0 \\
u(0,x)= u_0(x),\quad u_t(0,x)=u_1(x), &\quad x\in \R^n \\
v(0,x)= v_0(x),\quad v_t(0,x)=v_1(x), &\quad x\in \R^n.
\end{cases}
\end{equation}
Namely, the Cauchy problem for weakly coupled systems of semi-linear structurally damped $\sigma$-evolution equations
\begin{equation}\label{StructuralSystem_sigma}
\begin{cases}
u_{tt}+ (-\Delta)^\sigma u+ (-\Delta)^\delta u_t=|v|^p, &\quad x\in \R^n,\, t \ge 0 \\
v_{tt}+ (-\Delta)^\sigma v+ (-\Delta)^\delta v_t=|u|^p, &\quad x\in \R^n,\, t \ge 0 \\
u(0,x)= u_0(x),\quad u_t(0,x)=u_1(x), &\quad x\in \R^n \\
v(0,x)= v_0(x),\quad v_t(0,x)=v_1(x), &\quad x\in \R^n
\end{cases}
\end{equation}
with $\sigma\ge 1$ and $\delta \in [0,\sigma]$ was investigated in \cite{Dao2022}. The main aim of the cited paper is to prove global (in time) existence of small data Sobolev solutions to (\ref{StructuralSystem_sigma}) by assuming additional $L^m$ regularity on the initial data, with $m \in [1,2)$, together with using $(L^m \cap L^2)- L^2$ and $L^2- L^2$ estimates for solutions to the corresponding linear Cauchy problems. In addition, both the nonexistence of global solutions and lifespan estimates have been discussed in \cite{Dao2022} if $\sigma$ and $\delta$ are integers. We can say that the system (\ref{System_sigma}) is actually a special case of the problem studied in the very recent manuscript \cite{MezadekMezadekReissig2023}. By developing some $L^m-L^q$ estimates for Sobolev solutions to (\ref{LinearEq}), the authors in \cite{MezadekMezadekReissig2023} have verified the admissile range of power exponents $p,q$ of the nonlinear terms that allow small data Sobolev solutions to (\ref{System_sigma}) to exist in the whole space. Confining the data to the energy space and attaching appropriate high regularity to it, they compared the admissible range of $p,q$ with the modified Fujita exponent to understand how to relax the desired exponents. By choosing suitable integrability and regularity for the initial data, they have explored the blow-up results for (\ref{System_sigma}) and stated some estimates for lifespan. However, one recognizes that the achieved lifespan estimates in \cite{MezadekMezadekReissig2023} seem to be not sharp, i.e. there exists a gap between its upper bound and lower bound.

To the best of author's knowledge, there is no study on the multi-component weakly coupled systems (\ref{MainSytem}), i.e. $k>2$ in the literature which published before. Motivated strongly by the paper \cite{NishiharaWakasugi2015} in there taking account of the weakly coupled system of multi-semi-linear damped wave equations, the main goals of this paper are as follows: Under the supercritical condition
$$ \max\{\gamma_1,\gamma_2,\cdots,\gamma_{\mathtt{k}}\}< \frac{n}{2\sigma}, $$
we would like to prove the global solution existence to (\ref{MainSytem}) by assuming suitably small initial data and using some obtained estimates in \cite{DaoMichihisa2020}. We want to underline that our proof of global solutions is not a simple generalization of the proof in previous studies. More precisely, the core idea of this paper is that the calculation process of allowing multi-loss of decay plays an essential role in the existence of global solutions, which never appears in \cite{NishiharaWakasugi2015}. When the subcritical condition
$$ \max\{\gamma_1,\gamma_2,\cdots,\gamma_{\mathtt{k}}\}> \frac{n}{2\sigma} $$
holds, combined with appropriate conditions for the initial data and based on a contradiction argument, the blow-up phenomenon of solutions to (\ref{MainSytem}) occurs. We inherit the method of dealing with fractional Laplace operator in paper \cite{DaoReissig2021}, then apply it to derive the blow-up result and the upper bound estimates for lifespan simultaneously. Among other things, the next novelty is that we have succeeded in getting its lower bound to establish sharp estimates for lifespan of solutions to (\ref{MainSytem}) by developing the method of used in the paper \cite{ChenDao2023}. This means that our results come to fill out the gap related to lifespan estimates appearing in \cite{MezadekMezadekReissig2023}. Furthermore, our stratergy also works well for the problem in \cite{NishiharaWakasugi2015} to conclude a desired sharp result for lower bound of lifespan, which was a lack of \cite{NishiharaWakasugi2015}. Finally, we leave open the system (\ref{MainSytem}) in the critical condition
$$ \max\{\gamma_1,\gamma_2,\cdots,\gamma_{\mathtt{k}}\}= \frac{n}{2\sigma} $$
by the fact that to our knowledge there is no perfect way to deal with this situation so far. In other words, the main difficulty causes from the assumption of any fractional number $\sigma$ and taking account of the multi-component system.\medskip

\noindent\textbf{Notations:} Before stating our main results, we give the following notations which are used throughout this paper.
\begin{itemize}[leftmargin=*]
\item We write $f\lesssim g$ when there exists a constant $C>0$ such that $f\le Cg$, and $f \approx g$ when $g\lesssim f\lesssim g$.
\item As usual, the spaces $H^s$ and $\dot{H}^s$, with $s \ge 0$, stand for Bessel and Riesz potential spaces based on $L^2$ spaces. Here $\big<D\big>^s$ and $|D|^s$ denote the pseudo-differential operators with symbols $\big<\xi\big>^s$ and $|\xi|^s$, respectively.
\item For a given number $s \in \R$, we denote by
$$ [s]:= \max \big\{k \in \Z \,\, : \,\, k\le s \big\} \quad \text{ and }\quad [s]^+:= \max\{s,0\}, $$
its integer part and its positive part, respectively.
\item We introduce the space
$\mathcal{D}:= \big(L^1 \cap H^\sigma\big) \times \big(L^1 \cap L^2\big)$ with the norm
$$\|(\phi_0,\phi_1)\|_{\mathcal{D}}:=\|\phi_0\|_{L^1}+ \|\phi_0\|_{H^\sigma}+ \|\phi_1\|_{L^1}+ \|\phi_1\|_{L^2}, \quad \text{ where }\sigma\ge 1. $$
\item Let us denote the matrix
$$ P:=
\begin{pmatrix}
0&0&\cdots&0&p_1 \\
p_2&0&\cdots&0&0 \\
0&p_3&\cdots&0&0 \\
\vdots &\vdots& &\vdots &\vdots \\
0&0&\cdots&p_{\mathtt{k}}&0 \\
\end{pmatrix} $$
and the column vector $\gamma= (\gamma_1,\gamma_2,\cdots,\gamma_{\mathtt{k}})^{\mathtt{t}}:=(P-I_{\mathtt{k}})^{-1}(\underbrace{1,1,\cdots,1}_{\mathtt{k} \text{ times}})^{\mathtt{t}}$, where $I_{\mathtt{k}}$ and $(\alpha_1,\alpha_2,\cdots,\alpha_{\mathtt{k}})^{\mathtt{t}}$ stand for the identity matrix and the transposition vector of $(\alpha_1,\alpha_2,\cdots,\alpha_{\mathtt{k}})$, respectively.
\end{itemize}

\begin{remark}
\fontshape{n}
\selectfont
We want to point out that due to the assumption $p_\ell>1$ with $\ell=1,2,\cdots,\mathtt{k}$, it is obvious to recognize that
$$ \text{det}(P-I_{\mathtt{k}})= (-1)^{\mathtt{k}+1}\big(p_1p_2 \cdots p_\mathtt{k}-1\big) \neq 0. $$
So, the inverse matrix $(P-I_{\mathtt{k}})^{-1}$ exists. This means that the above vector $\gamma$ is well-defined. Before stating our main results, without loss of generality we assume that
$$ \gamma_{\mathtt{k}}= \max\{\gamma_1,\gamma_2,\cdots,\gamma_{\mathtt{k}}\}. $$
A straightforward calculation gives
$$ \gamma_{\mathtt{k}}= \frac{1+ p_{\mathtt{k}}+ p_{\mathtt{k}-1}p_{\mathtt{k}}+\cdots+ p_2p_3\cdots p_{\mathtt{k}}}{p_1p_2 \cdots p_\mathtt{k}-1}. $$
\end{remark}

Our main results read as follows.

\begin{theorem}[\textbf{Global existence}] \label{Global-existence.Thoerem}
Let $\sigma \ge 1$ and $1\le n\le 2\sigma$. We assume that $p_\ell\ge 2$ for any $\ell=1,2,\cdots,\mathtt{k}$ and satisfy the following conditions:
\begin{equation} \label{Cond.p_1}
p_1\le 1+\f{2\sigma}{n},
\end{equation}
\begin{equation} \label{Cond.p_l}
\f{p_1p_2 \cdots p_\ell-1}{1+ p_\ell+ p_{\ell-1}p_\ell+\cdots+ p_2p_3\cdots p_\ell}\le \f{2\sigma}{n}\quad \text{ with }\ell=2,3,\cdots,\mathtt{k}-1.
\end{equation}
In addition, we suppose the following condition:
\begin{equation} \label{exponent1A1}
\gamma_{\mathtt{k}}< \frac{n}{2\sigma}.
\end{equation}
Then, there exists a constant $\e_0>0$ such that for any small data
$$ \big((u_{01},u_{11}),(u_{02},u_{12}),\cdots,(u_{0\mathtt{k}},u_{1\mathtt{k}})\big) \in \underbrace{\mathcal{D} \times \mathcal{D} \times \cdots \times \mathcal{D}}_{\mathtt{k} \text{ times}} $$
satisfying the assumption
$$ \mathcal{I}[(u_{01},u_{11}),(u_{02},u_{12}),\cdots,(u_{0\mathtt{k}},u_{1\mathtt{k}})]:= \sum_{\ell=1}^{\mathtt{k}}\|(u_{0\ell},u_{1\ell})\|_{\mathcal{D}} \le \e_0, $$
the Cauchy problem \eqref{MainSytem} admits a unique global (in time) small data Sobolev solution
$$ (u_1,u_2,\cdots,u_{\mathtt{k}}) \in \Big(\mathcal{C}\big([0,\ity),H^\sigma\big)\Big)^{\mathtt{k}}. $$
Moreover, the following estimates hold for any $\ell=1,2,\cdots,\mathtt{k}$:
\begin{align*}
\|u_\ell(t,\cdot)\|_{L^2} &\lesssim (1+t)^{-\frac{n}{4\sigma}+\e_\ell} \mathcal{I}[(u_{01},u_{11}),(u_{02},u_{12}),\cdots,(u_{0\mathtt{k}},u_{1\mathtt{k}})], \\
\big\||D|^{\sigma} u_\ell(t,\cdot)\big\|_{L^2} &\lesssim (1+t)^{-\frac{n}{4\sigma}-\frac{1}{2}+\e_\ell} \mathcal{I}[(u_{01},u_{11}),(u_{02},u_{12}),\cdots,(u_{0\mathtt{k}},u_{1\mathtt{k}})],
\end{align*}
where the sequence $\{\e_\ell\}_{\ell=1}^{\mathtt{k}}$ is defined by
\begin{equation} \label{epsilon.formula-1}
\begin{cases}
\e_1= 1-\frac{n}{2\sigma}(p_1-1)+\e, \\
\e_\ell= 1-\frac{n}{2\sigma}(p_\ell-1)+p_\ell\e_{\ell-1}\quad \text{ for }\ell=2,3,\cdots,\mathtt{k}-1, \\
\e_{\mathtt{k}}=0,
\end{cases}
\end{equation}
for any arbitrarily small constant $\e>0$.
\end{theorem}

\begin{remark}
\fontshape{n}
\selectfont
From (\ref{epsilon.formula-1}), using an induction argument we compute
\begin{equation} \label{epsilon.formula-2}
\begin{cases}
\e_1= 1-\frac{n}{2\sigma}(p_1-1)+\e , \\
\e_2= 1+p_2-\frac{n}{2\sigma}(p_1p_2-1)+p_2\e , \\
\,\,\vdots \\
\e_{\mathtt{k}-1}= 1+ p_{\mathtt{k}-1}+ p_{\mathtt{k}-2}p_{\mathtt{k}-1}+\cdots+ p_2p_3\cdots p_{\mathtt{k}-1}- \frac{n}{2\sigma}\big(p_1p_2\cdots p_{\mathtt{k}-1}\big)+p_2p_3\cdots p_{\mathtt{k}-1}\e.
\end{cases}
\end{equation}
Thanks to the conditions \eqref{Cond.p_1} and \eqref{Cond.p_l}, one recognizes that $\e_\ell>0$ for all $\ell=1,2,\cdots,\mathtt{k}-1$, which represent some loss of decay of solutions to \eqref{MainSytem} in comparison with the corresponding linear equation \eqref{LinearEq}.
\end{remark}

\begin{theorem}[\textbf{Blow-up}] \label{Blow-up.Thoerem}
Let $\sigma \ge 1$. Assume that the initial data belong to the class
$$ \big(I+(-\Delta)^\sigma\big) u_{0\ell} \in L^1 \quad \text{ and }\quad u_{1\ell} \in L^1, $$
where $I$ represents the identity operator, and satisfy the relations
\begin{equation} \label{blow-up.1}
\int_{\R^n} \big(I+(-\Delta)^\sigma\big) u_{0\ell}(x) \,dx> 0 \quad \text{ and }\quad \int_{\R^n} u_{1\ell}(x) \,dx>0
\end{equation}
with $\ell=1,2,\cdots,\mathtt{k}$. Then, there is no global (in time) Sobolev solution
$$ (u_1,u_2,\cdots,u_{\mathtt{k}}) \in \big(\mathcal{C}\big([0,\ity),L^2\big)\big)^{\mathtt{k}} $$
to \eqref{MainSytem} if the following condition holds:
\begin{equation} \label{blow-up.2}
\gamma_{\mathtt{k}}> \frac{n}{2\sigma}.
\end{equation}
\end{theorem}

\begin{remark}
\fontshape{n}
\selectfont
From the conditions \eqref{exponent1A1} in Theorem \ref{Global-existence.Thoerem} and \eqref{blow-up.2} in Theorem \ref{Blow-up.Thoerem}, we understand that the condition
$$ \gamma_{\mathtt{k}}= \frac{n}{2\sigma}, \text{ in general, } \max\{\gamma_1,\gamma_2,\cdots,\gamma_{\mathtt{k}}\}= \frac{n}{2\sigma} $$
play an important role as a threshold between the global (in time) existence of small data solutions and the blow-up behavior of solutions in finite time. Hence, we would say that the above condition is really the critical condition for power multi-exponents in \eqref{MainSytem}.
\end{remark}

\begin{theorem}[\textbf{Upper bound of lifespan}] \label{Lifespan.Thoerem-1}
Assume that we have the hypothesis as in Theorem \ref{Blow-up.Thoerem} together with the condition \eqref{blow-up.2} for the power exponents. Then, there exists a positive constant $\varepsilon_0>0$ such that for any $\varepsilon \in (0,\varepsilon_0]$ the following upper bound estimate for the lifespan of solutions to \eqref{MainSytem} holds:
\begin{equation} \label{Upper_Lifespan}
{\rm LifeSpan}(u_1,u_2,\cdots,u_{\mathtt{k}})  \le C\varepsilon^{-\frac{1}{\gamma_{\mathtt{k}}-\frac{n}{2\sigma}}},
\end{equation}
where $C$ is a positive constant independent of $\varepsilon$.
\end{theorem}

\begin{theorem}[\textbf{Lower bound of lifespan}]  \label{Lifespan.Thoerem-2}
Let $\sigma \ge 1$. We suppose the initial data are taken from
$$ \big((u_{01},u_{11}),(u_{02},u_{12}),\cdots,(u_{0\mathtt{k}},u_{1\mathtt{k}})\big) \in \underbrace{\mathcal{D} \times \mathcal{D} \times \cdots \times \mathcal{D}}_{\mathtt{k} \text{ times}}. $$
Moreover, the power exponents $p_\ell$ with $\ell=1,\cdots, \mathtt{k}$ enjoy the condition \eqref{blow-up.2} together with the assumptions 
\begin{align}
&2 \le p_1,\,p_2,\cdots,p_{\mathtt{k}}< \ity &\quad &\text{ if }\, n\le 2\sigma, \label{GN1A1} \\
&2 \le p_1,\,p_2,\cdots,p_{\mathtt{k}}\le \frac{n}{n-2\sigma} &\quad &\text{ if }\, 2\sigma< n\le 4\sigma. \label{GN1A2}
\end{align}
Then, there exists a positive constant $\varepsilon_0$ such that for any $\e \in (0,\e_0]$ the following lower bound estimate for the lifespan of solutions to \eqref{MainSytem} holds:
\begin{equation} \label{Lower_Lifespan}
{\rm LifeSpan}(u_1,u_2,\cdots,u_{\mathtt{k}}) \ge c\varepsilon^{-\frac{1}{\gamma_{\mathtt{k}}-\frac{n}{2\sigma}}},
\end{equation}
where $c=c(n,u_{01},u_{11},\cdots,u_{0\mathtt{k}},u_{1\mathtt{k}})$ is a positive constant independent of $\varepsilon$.
\end{theorem}

\begin{remark}
\fontshape{n}
\selectfont
Connecting the two obtained results \eqref{Upper_Lifespan} and \eqref{Lower_Lifespan} in Theorems \ref{Lifespan.Thoerem-1} and \ref{Lifespan.Thoerem-2}, respectively, we conclude that the sharp lifespan estimates for solutions to the system \eqref{MainSytem} in the subcritical case \eqref{blow-up.2} are defined by 
$$ {\rm LifeSpan}(u_1,u_2,\cdots,u_{\mathtt{k}}) \sim \varepsilon^{-\frac{1}{\gamma_{\mathtt{k}}-\frac{n}{2\sigma}}}, $$
in general,
$$ {\rm LifeSpan}(u_1,u_2,\cdots,u_{\mathtt{k}}) \sim \varepsilon^{-\frac{1}{\max\{\gamma_1,\gamma_2,\cdots,\gamma_{\mathtt{k}}\}-\frac{n}{2\sigma}}}. $$
\end{remark}

%===========================================
\section{Global existence}\label{Sec.2}
\subsection{Useful estimates}\label{Sec.2.1}
In order to show the proof of Theorem \ref{Global-existence.Thoerem}, let us prepare several useful ingredients in the following propositions.
\begin{proposition}[\textbf{Linear estimates}, \cite{DaoMichihisa2020}] \label{prop.1.4.1}
Let $r \in [1,2)$. Then, the Sobolev solutions to (\ref{LinearEq}) satisfy the $(L^r \cap L^2)-L^2$ estimates
\begin{align*}
\big\|\partial_t^j |D|^s u_\ell(t,\cdot)\big\|_{L^2} &\lesssim (1+t)^{-\frac{n}{2\sigma}(\frac{1}{r}-\frac{1}{2})- \frac{s}{2\sigma}-j} \big(\|u_{0\ell}\|_{L^r \cap H^s}+ \|u_{1\ell}\|_{L^r \cap H^{[s+2(j-1)\sigma]^+}}\big),
\end{align*}
and the $L^2-L^2$ estimates
$$ \big\|\partial_t^j |D|^s u_\ell(t,\cdot)\big\|_{L^2} \lesssim (1+t)^{-\frac{s}{2\sigma}-j}\big(\|u_{0\ell}\|_{H^s}+ \|u_{1\ell}\|_{H^{[s+2(j-1)\sigma]^+}}\big), $$
for any $s\ge 0$, $j=0,1$, $\ell= 1,2,\cdots,\mathtt{k}$ and for all space dimensions $n\ge 1$.
\end{proposition}

\begin{proposition}[\textbf{Fractional Gagliardo-Nirenberg inequality}, \cite{ReissigEbert2018}] \label{fractionalGagliardoNirenberg}
Let $1<q,\,q_1,\,q_2<\infty$, $s>0$ and $a\in [0,s)$. Then, it holds for all $u\in L^{q_1} \cap \dot{H}^s_{q_2}$:
$$ \|u\|_{\dot{H}^a_q}\lesssim \|u\|_{L^{q_1}}^{1-\theta}\, \|u\|_{\dot{H}^s_{q_2}}^\theta, $$
where
$$ \theta=\theta(q,q_1,q_2,a,s,n)= \displaystyle\frac{\frac{1}{q_1}-\frac{1}{q}+\frac{a}{n}}{\frac{1}{q_1}-\frac{1}{q_2}+\frac{s}{n}} \,\,\text{ together with the condition }\,\, \f{a}{s}\leq \theta\leq 1. $$
\end{proposition}

%===========================================
\subsection{Proof of Theorem \ref{Global-existence.Thoerem}}
First of all, we introduce the solution space
$$ X(t):= \Big(\mathcal{C}\big([0,t],H^\sigma\big)\Big)^{\mathtt{k}} $$
with the norm
\begin{align*}
\|u\|_{X(t)}&= \|(u_1,u_2,\cdots,u_{\mathtt{k}})\|_{X(t)} \\
&:= \sup_{0\le \tau \le t} \Bigg(\sum_{\ell=1}^{\mathtt{k}-1}(1+\tau)^{\frac{n}{4\sigma}-\e_\ell}\|u_\ell(\tau,\cdot)\|_{L^2}+ \sum_{\ell=1}^{\mathtt{k}-1}(1+\tau)^{\frac{n}{4\sigma}+\frac{1}{2}-\e_\ell}\big\||D|^{\sigma} u_\ell(\tau,\cdot)\big\|_{L^2} \\
&\hspace{2cm}+ (1+\tau)^{\frac{n}{4\sigma}}\|u_{\mathtt{k}}(\tau,\cdot)\|_{L^2}+ (1+\tau)^{\frac{n}{4\sigma}+\frac{1}{2}}\big\||D|^{\sigma} u_{\mathtt{k}}(\tau,\cdot)\big\|_{L^2}\Bigg),
\end{align*}
where the sequence $\{\e_\ell\}_{\ell=1}^{\mathtt{k}-1}$ is defined as in \eqref{epsilon.formula-1}. Denoting by $K_0(t,x)$ and $K_1(t,x)$ the fundamental solutions to (\ref{LinearEq}), we can write the solutions of the corresponding linear Cauchy problems with vanishing right-hand sides to (\ref{MainSytem}) in the following form:
\begin{equation} \label{Sol.Rep-1}
u_\ell^{\rm ln}(t,x)=K_0(t,x)*_x u_{0\ell}(x)+K_1(t,x)*_x u_{1\ell}(x)
\end{equation}
for $\ell= 1,2,\cdots,\mathtt{k}$, where $*_x$ stands for the convolution with respect to the spatial variables $x$. Using Duhamel's principle we get the formal implicit representation of solutions to (\ref{MainSytem}) in the following form:
\begin{equation} \label{Sol.Rep-2}
\begin{cases}
u_1(t,x)=u_1^{\rm ln}(t,x)+ \dps\int_0^t K_1(t-\tau,x)*_x |u_\mathtt{k}(\tau,x)|^{p_1}d\tau =u_1^{\rm ln}(t,x)+u_1^{\rm nl}(t,x) \\
u_\ell(t,x)=u_\ell^{\rm ln}(t,x)+ \dps\int_0^t K_1(t-\tau,x)*_x |u_{\ell-1}(\tau,x)|^{p_\ell}d\tau =u_\ell^{\rm ln}(t,x)+u_\ell^{\rm nl}(t,x)
\end{cases}
\end{equation}
for $\ell= 2,3,\cdots,\mathtt{k}$. We define the following operator for all $t>0$:
\begin{align*}
N:X(t) &\rightarrow X(t) \\
N[u](t,x) &=u^{\rm ln}(t,x)+u^{\rm nl}(t,x).
\end{align*}
We will show that the operator $N$ fulfills the following two inequalities:
\begin{align}
\|N[u]\|_{X(t)} &\lesssim \sum_{\ell=1}^{\mathtt{k}}\|(u_{0\ell},u_{1\ell})\|_{\mathcal{D}}+ \sum_{\ell=1}^{\mathtt{k}} \|u\|^{p_{\ell}}_{X(t)}, \label{operator.1} \\
\|N[u]-N[v]\|_{X(t)} &\lesssim \|u-v\|_{X(t)}\left(\sum_{\ell=1}^{\mathtt{k}} \|u\|^{p_{\ell}-1}_{X(t)}+\sum_{\ell=1}^{\mathtt{k}} \|v\|^{p_{\ell}-1}_{X(t)}\right). \label{operator.2}
\end{align}
Then, we can conclude global existence results of small data solution by applying Banach's fixed point theorem. 

\textit{At first, let us verify the inequality \eqref{operator.1}}. Obviously, from Proposition \ref{prop.1.4.1} and the definition of the norm in $X(t)$ we get
$$ \|u^{\rm ln}\|_{X(t)}\lesssim \sum_{\ell=1}^{\mathtt{k}}\|(u_{0\ell},u_{1\ell})\|_{\mathcal{D}}. $$
Then, it is reasonable to prove the following inequality instead of (\ref{operator.1}):
\begin{equation} \label{operator.3}
 \|u^{\rm nl}\|_{X(t)}\lesssim \sum_{\ell=1}^{\mathtt{k}}\|u\|^{p_{\ell}}_{X(t)}.
 \end{equation}
Indeed, at the first stage we employ the fractional Gagliardo-Nirenberg inequality from Proposition \ref{fractionalGagliardoNirenberg} to gain the following estimates:
\begin{align*}
\big\||u_\ell(\tau,\cdot)|^{p_{\ell+1}}\big\|_{L^1} &= \big\||u_\ell(\tau,\cdot)\|^{p_{\ell+1}}_{L^{p_{\ell+1}}} \\
&\lesssim \Big(\|u_\ell(\tau,\cdot)\|^{1-\theta_{1\ell}}_{L^2}\, \big\||D|^{\sigma} u_\ell(\tau,\cdot)\big\|^{\theta_{1\ell}}_{L^2}\Big) \qquad \text{ with } \theta_{1\ell}= \frac{n}{\sigma}\left(\frac{1}{2}-\frac{1}{p_{\ell+1}}\right) \in [0,1] \\
&\lesssim \Big((1+\tau)^{-\frac{n}{4\sigma}-\frac{1}{2}\theta_{1\ell}+\e_\ell}\Big)^{p_{\ell+1}}\|u\|^{p_{\ell+1}}_{X(t)} \\
&= (1+\tau)^{-\frac{n}{2\sigma}(p_{\ell+1}-1)+p_{\ell+1}\e_\ell} \|u\|^{p_{\ell+1}}_{X(t)},
\end{align*}
and
\begin{align*}
\big\||u_\ell(\tau,\cdot)|^{p_{\ell+1}}\big\|_{L^2} &= \big\||u_\ell(\tau,\cdot)\|^{p_{\ell+1}}_{L^{2p_{\ell+1}}} \\
&\lesssim \Big(\|u_\ell(\tau,\cdot)\|^{1-\theta_{2\ell}}_{L^2}\, \big\||D|^{\sigma} u_\ell(\tau,\cdot)\big\|^{\theta_{2\ell}}_{L^2}\Big) \qquad \text{ with } \theta_{2\ell}= \frac{n}{\sigma}\left(\frac{1}{2}-\frac{1}{2p_{\ell+1}}\right) \in [0,1] \\
&\lesssim \Big((1+\tau)^{-\frac{n}{4\sigma}-\frac{1}{2}\theta_{2\ell}+\e_\ell}\Big)^{p_{\ell+1}}\|u\|^{p_{\ell+1}}_{X(t)} \\
&= (1+\tau)^{-\frac{n}{2\sigma}(p_{\ell+1}-\frac{1}{2})+p_{\ell+1}\e_\ell} \|u\|^{p_{\ell+1}}_{X(t)},
\end{align*}
which implies immediately
\begin{align*}
\big\||u_\ell(\tau,\cdot)|^{p_{\ell+1}}\big\|_{L^1\cap L^2} &= \big\||u_\ell(\tau,\cdot)|^{p_{\ell+1}}\big\|_{L^1}+ \big\||u_\ell(\tau,\cdot)|^{p_{\ell+1}}\big\|_{L^2} \\ 
&\lesssim (1+\tau)^{-\frac{n}{2\sigma}(p_{\ell+1}-1)+p_{\ell+1}\e_\ell} \|u\|^{p_{\ell+1}}_{X(t)}
\end{align*}
for $\ell= 1,2,\cdots,\mathtt{k}-1$. Simultaneously, one also derives
\begin{align*}
\big\||u_{\mathtt{k}}(\tau,\cdot)|^{p_1}\big\|_{L^1} &\lesssim \Big(\|u_{\mathtt{k}}(\tau,\cdot)\|^{1-\theta_{1\ell}}_{L^2}\, \big\||D|^{\sigma} u_{\mathtt{k}}(\tau,\cdot)\big\|^{\theta_{1\mathtt{k}}}_{L^2}\Big) \lesssim (1+\tau)^{-\frac{n}{2\sigma}(p_1-1)} \|u\|^{p_1}_{X(t)}, \\
\big\||u_{\mathtt{k}}(\tau,\cdot)|^{p_1}\big\|_{L^2} &\lesssim \Big(\|u_{\mathtt{k}}(\tau,\cdot)\|^{1-\theta_{2\mathtt{k}}}_{L^2}\, \big\||D|^{\sigma} u_{\mathtt{k}}(\tau,\cdot)\big\|^{\theta_{2\mathtt{k}}}_{L^2}\Big) \lesssim (1+\tau)^{-\frac{n}{2\sigma}(p_1-\frac{1}{2})} \|u\|^{p_1}_{X(t)},
\end{align*}
where
$$ \theta_{1\mathtt{k}}= \frac{n}{\sigma}\left(\frac{1}{2}-\frac{1}{p_1}\right) \in [0,1]\quad \text{ and }\quad \theta_{2\mathtt{k}}= \frac{n}{\sigma}\left(\frac{1}{2}-\frac{1}{2p_1}\right) \in [0,1]. $$
Thus, it follows that
\begin{align*}
\big\||u_{\mathtt{k}}(\tau,\cdot)|^{p_1}\big\|_{L^1\cap L^2} &= \big\||u_{\mathtt{k}}(\tau,\cdot)|^{p_1}\big\|_{L^1}+ \big\||u_{\mathtt{k}}(\tau,\cdot)|^{p_1}\big\|_{L^2}\\ 
&\lesssim (1+\tau)^{-\frac{n}{2\sigma}(p_1-1)} \|u\|^{p_1}_{X(t)}.
\end{align*}
With $\ell= 1,2,\cdots,\mathtt{k}$, it is clear to check that all these conditions for $\theta_{1\ell}$ and $\theta_{2\ell}$ hold as long as the assumptions $n\le 2\sigma$ and $p_\ell\ge 2$ are provided. Let us now come back to estimate some nonlinear terms in the second stage. The main stratergy is to use the $(L^1 \cap L^2)- L^2$ estimates if $\tau \in [0,t/2]$ and the $L^2-L^2$ estimates if $\tau \in [t/2,t]$ for solutions to the corresponding linear equation (\ref{LinearEq}) from Proposition \ref{prop.1.4.1} in the following ways:
\begin{align*}
\big\||D|^{\mathtt{m}\sigma} u^{\text{nl}}_1(\tau,\cdot)\big\|_{L^2} &\lesssim \int_0^{t/2}(1+t-\tau)^{-\frac{n}{4\sigma}-\frac{\mathtt{m}}{2}}\big\||u_\mathtt{k}(\tau,\cdot)|^{p_1}\big\|_{L^1\cap L^2}d\tau \\
&\qquad + \int_{t/2}^t (1+t-\tau)^{-\frac{\mathtt{m}}{2}}\big\||u_\mathtt{k}(\tau,\cdot)|^{p_1}\big\|_{L^2}d\tau \\
&\lesssim (1+t)^{-\frac{n}{4\sigma}-\frac{\mathtt{m}}{2}}\|u\|^{p_1}_{X(t)} \int_0^{t/2}(1+\tau)^{-\frac{n}{2\sigma}(p_1-1)}d\tau \\
&\qquad + (1+t)^{-\frac{n}{2\sigma}(p_1-\frac{1}{2})} \|u\|^{p_1}_{X(t)} \int_{t/2}^t (1+t-\tau)^{-\frac{\mathtt{m}}{2}}d\tau \\
&\lesssim (1+t)^{-\frac{n}{4\sigma}-\frac{\mathtt{m}}{2}+\e_1}\|u\|^{p_1}_{X(t)}
\end{align*}
with $\mathtt{m}=0,1$, where we notice that the condition \eqref{Cond.p_1} leads to
$$ -\frac{n}{2\sigma}(p_1-1)\ge -1. $$
Analogously, we arrive at
\begin{align*}
\big\||D|^{\mathtt{m}\sigma} u^{\text{nl}}_\ell(\tau,\cdot)\big\|_{L^2} &\lesssim \int_0^{t/2}(1+t-\tau)^{-\frac{n}{4\sigma}-\frac{\mathtt{m}}{2}}\big\||u_{\ell-1}(\tau,\cdot)|^{p_\ell}\big\|_{L^1\cap L^2}d\tau \\
&\qquad + \int_{t/2}^t (1+t-\tau)^{-\frac{\mathtt{m}}{2}}\big\||u_{\ell-1}(\tau,\cdot)|^{p_\ell}\big\|_{L^2}d\tau \\
&\lesssim (1+t)^{-\frac{n}{4\sigma}-\frac{\mathtt{m}}{2}}\|u\|^{p_\ell}_{X(t)} \int_0^{t/2}(1+\tau)^{-\frac{n}{2\sigma}(p_\ell-1)+p_\ell\e_{\ell-1}}d\tau \\
&\qquad + (1+t)^{-\frac{n}{2\sigma}(p_\ell-\frac{1}{2})+p_\ell\e_{\ell-1}} \|u\|^{p_\ell}_{X(t)} \int_{t/2}^t (1+t-\tau)^{-\frac{\mathtt{m}}{2}}d\tau \\
&\lesssim (1+t)^{-\frac{n}{4\sigma}-\frac{\mathtt{m}}{2}+\e_\ell}\|u\|^{p_\ell}_{X(t)}
\end{align*}
with $\mathtt{m}=0,1$ and $\ell= 2,\cdots,\mathtt{k}-1$. Here observing from \eqref{epsilon.formula-1} and \eqref{epsilon.formula-2} that
\begin{align*}
\e_{\ell}&= 1-\frac{n}{2\sigma}(p_\ell-1)+p_\ell\e_{\ell-1} \\ 
&= 1+ p_{\ell}+ p_{\ell-1}p_{\ell}+\cdots+ p_2p_3\cdots p_{\ell}- \f{n}{2\sigma}\big(p_1p_2 \cdots p_\ell-1\big)+p_2p_3 \cdots p_\ell\e >0
\end{align*}
by the help of the condition \eqref{Cond.p_l}, we deduce immediately  the relation
$$ -\frac{n}{2\sigma}(p_\ell-1)+p_\ell\e_{\ell-1}> -1. $$
In the same manner, one also concludes that
\begin{align*}
\big\||D|^{\mathtt{m}\sigma} u^{\text{nl}}_\mathtt{k}(\tau,\cdot)\big\|_{L^2} &\lesssim \int_0^{t/2}(1+t-\tau)^{-\frac{n}{4\sigma}-\frac{\mathtt{m}}{2}}\big\||u_{\mathtt{k}-1}(\tau,\cdot)|^{p_\mathtt{k}}\big\|_{L^1\cap L^2}d\tau \\
&\qquad + \int_{t/2}^t (1+t-\tau)^{-\frac{\mathtt{m}}{2}}\big\||u_{\mathtt{k}-1}(\tau,\cdot)|^{p_\mathtt{k}}\big\|_{L^2}d\tau \\
&\lesssim (1+t)^{-\frac{n}{4\sigma}-\frac{\mathtt{m}}{2}}\|u\|^{p_\mathtt{k}}_{X(t)} \int_0^{t/2}(1+\tau)^{-\frac{n}{2\sigma}(p_\mathtt{k}-1)+p_\mathtt{k}\e_{\mathtt{k}-1}}d\tau \\
&\qquad + (1+t)^{-\frac{n}{2\sigma}(p_\mathtt{k}-\frac{1}{2})+p_\mathtt{k}\e_{\mathtt{k}-1}} \|u\|^{p_\mathtt{k}}_{X(t)} \int_{t/2}^t (1+t-\tau)^{-\frac{\mathtt{m}}{2}}d\tau \\
&\lesssim (1+t)^{-\frac{n}{4\sigma}-\frac{\mathtt{m}}{2}}\|u\|^{p_\mathtt{k}}_{X(t)},
\end{align*}
provided that the following condition is fulfilled:
\begin{equation} \label{epsilon.k}
-\frac{n}{2\sigma}(p_\mathtt{k}-1)+p_\mathtt{k}\e_{\mathtt{k}-1}< -1.
\end{equation}
Recalling \eqref{epsilon.formula-2} we get
\begin{equation} \label{epsilon.k-1}
\e_{\mathtt{k}-1}= 1+ p_{\mathtt{k}-1}+ p_{\mathtt{k}-2}p_{\mathtt{k}-1}+\cdots+ p_2p_3\cdots p_{\mathtt{k}-1}- \f{n}{2\sigma}\big(p_1p_2\cdots p_\mathtt{k}-1\big)+ p_2p_3\cdots p_\mathtt{k}\e.
\end{equation}
For this reason, plugging (\ref{epsilon.k-1}) into (\ref{epsilon.k}) we achieve the following condition:
$$ \frac{1+ p_{\mathtt{k}}+ p_{\mathtt{k}-1}p_{\mathtt{k}}+\cdots+ p_2p_3\cdots p_{\mathtt{k}}}{p_1p_2p_3\cdots p_{\mathtt{k}-1}}< \frac{n}{2\sigma}, \quad \text{ that is, }\quad \gamma_{\mathtt{k}}< \frac{n}{2\sigma}, $$
where we have utilized the assumption of arbitrarily small constant $\e>0$. From the definition of the norm in $X(t)$ we can conclude the inequality (\ref{operator.3}).

\textit{Now let us turn to prove the inequality \eqref{operator.2}.} First of all, taking into consideration two elements $u=(u_1,u_2,\cdots, u_\mathtt{k})$ and $v=(v_1,v_2,\cdots, v_\mathtt{k})$ from $X(t)$ we get
$$ N[u](t,x)- N[v](t,x)= \Big(u_1^{\text{nl}}(t,x)- v_1^{\text{nl}}(t,x), u_2^{\text{nl}}(t,x)- v_2^{\text{nl}}(t,x), \cdots, u_\mathtt{k}^{\text{nl}}(t,x)- v_\mathtt{k}^{\text{nl}}(t,x)\Big). $$
Again, employing $(L^1 \cap L^2)- L^2$ estimates if $\tau \in [0,t/2]$ and the $L^2-L^2$ estimates if $\tau \in [t/2,t]$ from Proposition \ref{prop.1.4.1} we obtain the following estimates with $\mathtt{m}=0,1$:
\begin{align*}
\big\||D|^{\mathtt{m}\sigma} (u^{\text{nl}}_\ell- v^{\text{nl}}_\ell)(\tau,\cdot)\big\|_{L^2} &\lesssim \int_0^{t/2}(1+t-\tau)^{-\frac{n}{2\sigma}(\frac{1}{r}-\frac{1}{2})-\frac{\mathtt{m}}{2}}\big\||u_{\ell-1}(\tau,\cdot)|^{p_\ell}-|v_{\ell-1}(\tau,\cdot)|^{p_\ell}\big\|_{L^1\cap L^2}d\tau \\
&\quad + \int_{t/2}^t (1+t-\tau)^{-\frac{\mathtt{m}}{2}}\big\||u_{\ell-1}(\tau,\cdot)|^{p_\ell}-|v_{\ell-1}(\tau,\cdot)|^{p_\ell}\big\|_{L^2}d\tau.
\end{align*}
After applying H\"{o}lder's inequality, one establishes
\begin{align*}
\big\||u_{\ell-1}(\tau,\cdot)|^{p_\ell}-|v_{\ell-1}(\tau,\cdot)|^{p_\ell}\big\|_{L^1}\lesssim \|u_{\ell-1}(\tau,\cdot)-v_{\ell-1}(\tau,\cdot)\|_{L^{p_\ell}}\Big(\|u_{\ell-1}(\tau,\cdot)\|^{p_\ell-1}_{L^{p_\ell}}+\|v_{\ell-1}(\tau,\cdot)\|^{p_\ell-1}_{L^{p_\ell}}\Big), \\
\big\||u_{\ell-1}(\tau,\cdot)|^{p_\ell}-|v_{\ell-1}(\tau,\cdot)|^{p_\ell}\big\|_{L^2}\lesssim \|u_{\ell-1}(\tau,\cdot)-v_{\ell-1}(\tau,\cdot)\|_{L^{2p_\ell}}\Big(\|u_{\ell-1}(\tau,\cdot)\|^{p_\ell-1}_{L^{2p_\ell}}+\|v_{\ell-1}(\tau,\cdot)\|^{p_\ell-1}_{L^{2p_\ell}}\Big).
\end{align*}
Similarly to the proof of (\ref{operator.3}), we use the fractional Gagliardo-Nirenberg inequality from Proposition \ref{fractionalGagliardoNirenberg} in terms of estimating
$$ \|u_{\ell-1}(\tau,\cdot)-v_{\ell-1}(\tau,\cdot)\|_{L^\eta},\quad \|u_{\ell-1}(\tau,\cdot)\|_{L^\eta},\quad \|v_{\ell-1}(\tau,\cdot)\|_{L^\eta} $$
with $\eta=p_\ell$ or $\eta=2p_\ell$ and then repeat some arguments as we did to indicate (\ref{operator.3}). In this way, we may conclude the inequality (\ref{operator.2}), hence, the proof of Theorem \ref{Global-existence.Thoerem} is completed.

%===========================================
\section{Blow-up result}
\subsection{Auxiliary tools}
Before presenting the proof of Theorem \ref{Blow-up.Thoerem} in this section, let us recall the following auxiliary lemmas which are some essential properties of a modified test function from the quite recent paper \cite{DaoReissig2021} of the second author to prove the blow-up result.
\begin{lemma} \label{lemma2.1}
Let $m \in \N$ and $s \in [0,1)$. Then, the following estimates hold for any $q>n$ and for all $x \in \R^n$:
\begin{align*}
 \left|(-\Delta)^{m+s} \langle x\rangle^{-q}\right| \lesssim
\begin{cases}
\langle x\rangle^{-n-2m} & \text{ if }\, s=0, \\
\langle x\rangle ^{-n-2s} & \text{ if }\, s \in (0,1).
\end{cases}
\end{align*}
\end{lemma}

\begin{lemma} \label{lemma2.2}
Let $\nu \geqslant 1$ be a fractional number. Let $\varphi:= \varphi(x)= \langle x\rangle^{-q}$ for some $q>0$. For any $R>0$, let $\varphi_R$ be a function defined by
$$ \varphi_R(x):= \varphi(x/R)\quad  \text{ for all } x \in \R^n. $$
Then, $(-\Delta)^\nu (\varphi_R)$ satisfies the following scaling properties for all $x \in \R^n$:
\begin{equation*}
(-\Delta)^\nu (\varphi_R)(x)= R^{-2\nu} \left((-\Delta)^\nu \varphi \right)(x/R).
\end{equation*}
\end{lemma}

\begin{lemma} \label{lemma2.3}
Let $s \in \R$. Let $\varphi_1=\varphi_1(x) \in H^s(\R^n)$ and $\varphi_2=\varphi_2(x) \in H^{-s}(\R^n)$. Then, the following relation holds:
$$ \int_{\R^n}\varphi_1(x)\,\varphi_2(x)dx= \int_{\R^n}\widehat{\varphi}_1(\xi)\,\widehat{\varphi}_2(\xi)d\xi, $$
where we denote by $\widehat{\varphi}_1(\xi)$ and $\widehat{\varphi}_2(\xi)$, the Fourier transform of functions $\varphi_1(x)$ and $\varphi_2(x)$, respectively.
\end{lemma}

%===========================================
\subsection{Proof of Theorem \ref{Blow-up.Thoerem}}
Firstly, for any fractional number $\sigma\ge 1$ let us denote the following constant:
$$ \bar{\sigma}:=
\begin{cases}
1 &\text{ if }\sigma \text{ is an integer number}, \\
\sigma- [\sigma] &\text{ if }\sigma \text{ is not an integer number}.
\end{cases}$$
Now we introduce the radial space-dependent test function $\psi=\psi(x)$ by
\begin{equation} \label{key.estimate}
\psi(x):=\langle x\rangle^{-n-2\bar{\sigma}}=(1+|x|^2)^{-n/2-\bar{\sigma}}.
\end{equation}
Then, the application of Lemma \ref{lemma2.1} leads to the following estimate for any fractional number $\sigma\ge 1$:
\begin{equation}\label{key_Laplace}
\big|(-\Delta)^{\sigma} \psi(x)\big| \lesssim \langle x\rangle^{-n-2\bar{\sigma}}.
\end{equation}
Moreover, we choose the time-dependent test function $\eta= \eta(t)$ satisfying
\begin{align}
&1.\quad \eta \in \mathcal{C}_0^\ity([0,\ity)) \text{ and }
\eta(t)=\begin{cases}
1 &\quad \text{ if }\ \ 0 \le t \le 1/2, \\
\text{decreasing } &\quad \text{ if }\ \ 1/2 < t < 1, \\
0 &\quad \text{ if }\ \ t \ge 1,
\end{cases} \nonumber \\
&2.\quad \eta^{-\frac{\lambda'}{\lambda}}(t)\big(|\eta'(t)|^{\lambda'}+|\eta''(t)|^{\lambda'}\big) \le C \quad \text{ for any } t \in [1/2,1], \label{eta.Cond}
\end{align}
where $\lambda'$ is the conjugate of $\lambda>1$ and $C$ is a suitable positive constant.
\par Let $R$ be a large parameter in $[0,\ity)$. We define the test function
$$ \Phi_R(t,x):= \eta_R(t) \psi_R(x), $$
where $\eta_R(t):= \eta\big(R^{-2\sigma}t\big)$ and $\psi_R(x):= \psi\big(R^{-1}x\big)$. We define the following functionals with $\ell=1,2,\cdots,\mathtt{k}-1$:
$$ F_R[u_\ell]:= \int_0^{\ity}\int_{\R^n} |u_\ell(t,x)|^{p_{\ell+1}} \Phi_R(t,x)\,dxdt= \int_0^{R^{2\sigma}}\int_{\R^n}|u_\ell(t,x)|^{p_{\ell+1}} \Phi_R(t,x)\,dxdt $$
and
$$ F_R[u_{\mathtt{k}}]:= \int_0^{\ity}\int_{\R^n} |u_{\mathtt{k}}(t,x)|^{p_1} \Phi_R(t,x)\,dxdt= \int_0^{R^{2\sigma}}\int_{\R^n}|u_{\mathtt{k}}(t,x)|^{p_1} \Phi_R(t,x)\,dxdt. $$
Let us assume that $(u_1,u_2,\cdots,u_{\mathtt{k}})= \big(u_1(t,x),u_2(t,x),\cdots,u_{\mathtt{k}}(t,x)\big)$ is a global (in time) Sobolev solution from $\big(\mathcal{C}\big([0,\ity),L^2\big)\big)^{\mathtt{k}}$ to (\ref{MainSytem}). After multiplying the first equation of (\ref{MainSytem}) by $\Phi_R=\Phi_R(t,x)$, we carry out partial integration to derive
\begin{align}
&F_R[u_{\mathtt{k}}] + \int_{\R^n} \big((-\Delta)^\sigma u_{01}(x)+ u_{01}(x)+ u_{11}(x)\big)\psi_R(x)\,dx \nonumber \\
&\quad= \int_{R^{2\sigma}/2}^{R^{2\sigma}}\int_{\R^n}\eta''_R(t) \psi_R(x)u_1(t,x) \,dxdt + \int_0^{\ity}\int_{\R^n} \eta_R(t) \psi_R(x)\, (-\Delta)^{\sigma} u_1(t,x)\,dxdt \nonumber \\
&\qquad- \int_{R^{2\sigma}/2}^{R^{2\sigma}}\int_{\R^n} \eta'_R(t) \psi_R(x)\,u_1(t,x)\,dxdt - \int_{R^{2\sigma}/2}^{R^{2\sigma}}\int_{\R^n} \eta'_R(t) \psi_R(x)\,(-\Delta)^{\sigma} u_1(t,x)\,dxdt \nonumber \\
&\quad =: F_{1,R}[u_{\mathtt{k}}]+ F_{2,R}[u_{\mathtt{k}}]- F_{3,R}[u_{\mathtt{k}}]- F_{4,R}[u_{\mathtt{k}}]. \label{blow-up.t.1}
\end{align}
Applying H\"{o}lder's inequality with $\frac{1}{p_2}+ \frac{1}{p'_2}= 1$ we may estimate as follows:
\begin{align*}
\big|F_{1,R}[u_{\mathtt{k}}]\big| &\le \int_{R^{2\sigma}/2}^{R^{2\sigma}}\int_{\R^n} |u_1(t,x)|\, \big|\eta''_R(t)\big| \psi_R(x) \, dxdt \\
&\lesssim \left(\int_{R^{2\sigma}/2}^{R^{2\sigma}}\int_{\R^n} \Big|u_1(t,x)\phi^{\frac{1}{p_2}}_R(t,x)\Big|^{p_2} \,dxdt \right)^{\frac{1}{p_2}} \left(\int_{R^{2\sigma}/2}^{R^{2\sigma}}\int_{\R^n} \Big|\phi^{-\frac{1}{p_2}}_R(t,x) \eta''_R(t) \psi_R(x)\Big|^{p'_2}\, dxdt \right)^{\frac{1}{p'_2}} \\
&\lesssim F_R[u_1]^{\frac{1}{p_2}}\, \left(\int_{R^{2\sigma}/2}^{R^{2\sigma}}\int_{\R^n} \eta_R^{-\frac{p'_2}{p_2}}(t) \big|\eta''_R(t)\big|^{p'_2} \psi_R(x)\, dxdt\right)^{\frac{1}{p'_2}}.
\end{align*}
By the change of variables $\tilde{t}:= R^{-2\sigma}t$ and $\tilde{x}:= R^{-1}x$, a straight-forward calculation gives
\begin{align}
\big|F_{1,R}[u_{\mathtt{k}}]\big| &\lesssim F_R[u_1]^{\frac{1}{p_2}}\, R^{-4\sigma+ \frac{n+2\sigma}{p'_2}}\left(\int_{\R^n} \langle \tilde{x} \rangle^{-n-2\bar{\sigma}}\, d\tilde{x}\right)^{\frac{1}{p'_2}} \nonumber \\
&\lesssim R^{-4\sigma+ \frac{n+2\sigma}{p'_2}}F_R[u_1]^{\frac{1}{p_2}}. \label{blow-up.t.2}
\end{align}
Here we have used $\eta''_R(t)= R^{-2\sigma}\eta''(\tilde{t})$ and the assumption (\ref{eta.Cond}). In the same way, one also derives
\begin{equation}
\big|F_{3,R}[u_{\mathtt{k}}]\big| \lesssim R^{-2\sigma+ \frac{n+2\sigma}{p'_2}}F_R[u_1]^{\frac{1}{p_2}}, \label{blow-up.t.3}
\end{equation}
where noticing that the relation $\eta''_R(t)= R^{-2\sigma}\eta''(\tilde{t})$ holds and the assumption (\ref{eta.Cond}) has been utilized again. Now let us turn to estimate $F_{2,R}[u_{\mathtt{k}}]$ and $F_{4,R}[u_{\mathtt{k}}]$. First, by using $\psi_R \in H^{2\sigma}$ and $u_1 \in \mathcal{C}\big([0,\infty),L^2\big)$ we apply Lemma \ref{lemma2.3} to conclude the following relation:
$$ \int_{\R^n} \psi_R(x)\, (-\Delta)^{\sigma} u_1(t,x)\,dx= \int_{\R^n}|\xi|^{2\sigma}\widehat{\psi}_R(\xi)\,\widehat{u}_1(t,\xi)\,d\xi= \int_{\R^n} u_1(t,x)\, (-\Delta)^{\sigma}\psi_R(x)\,dx. $$
Hence, we obtain
\begin{align*}
F_{2,R}[u_{\mathtt{k}}] &= \int_0^{\ity}\int_{\R^n} \eta_R(t) u_1(t,x)\, (-\Delta)^{\sigma}\psi_R(x) \,dxdt, \\
F_{4,R}[u_{\mathtt{k}}] &= \int_{R^{2\sigma}/2}^{R^{2\sigma}}\int_{\R^n} \eta'_R(t) u_1(t,x)\,(-\Delta)^{\sigma}\psi_R(x)\, dxdt.
\end{align*}
Emloying again H\"{o}lder's inequality with $\frac{1}{p_2}+ \frac{1}{p'_2}= 1$ as we estimated $F_{1,R}[u_{\mathtt{k}}]$ leads to
$$ \big|F_{2,R}[u_{\mathtt{k}}]\big|\lesssim F_R[u_1]^{\frac{1}{p_2}}\, \left(\int_0^{R^{2\sigma}}\int_{\R^n} \eta_R(t) \psi^{-\frac{p'_2}{p_2}}_R(x)\, \big|(-\Delta)^{\sigma}\psi_R(x)\big|^{p'_2} \, dxdt\right)^{\frac{1}{p'_2}} $$
and
$$ \big|F_{4,R}[u_{\mathtt{k}}]\big|\lesssim F_R[u_1]^{\frac{1}{p_2}}\, \left(\int_{R^{2\sigma}/2}^{R^{2\sigma}}\int_{\R^n} \eta^{-\frac{p'_2}{p_2}}_R(t) \big|\eta'_R(t)\big|^{p'_2} \psi^{-\frac{p'_2}{p_2}}_R(x)\, \big|(-\Delta)^{\sigma}\psi_R(x)\big|^{p'_2}\,dxdt\right)^{\frac{1}{p'_2}}. $$
In order to control the above two integrals, the key tools rely on Lemma \ref{lemma2.2} and the estimate \eqref{key_Laplace}. Namely, at first carrying out the change of variables $\tilde{t}:= R^{-2\sigma}t$ and $\tilde{x}:= R^{-1}x$ we arrive at
\begin{align}
\big|F_{2,R}[u_{\mathtt{k}}]\big| &\lesssim F_R[u_1]^{\frac{1}{p_2}}\, R^{-2\sigma+ \frac{n+2\sigma}{p'_2}}\left(\int_0^{1}\int_{\R^n} \eta(\tilde{t}) \psi^{-\frac{p'_2}{p_2}}(\tilde{x})\, \big|(-\Delta)^{\sigma}(\psi)(\tilde{x})\big|^{p'_2}\, d\tilde{x}d\tilde{t}\right)^{\frac{1}{p'_2}} \nonumber \\
&\hspace{7cm} (\text{by  Lemma  }\ref{lemma2.2}) \nonumber \\
&\lesssim F_R[u_1]^{\frac{1}{p_2}}\, R^{-2\sigma+ \frac{n+2\sigma}{p'_2}}\left(\int_{\R^n} \psi^{-\frac{p'_2}{p_2}}(\tilde{x})\, \big|(-\Delta)^{\sigma}(\psi)(\tilde{x})\big|^{p'_2}\, d\tilde{x}\right)^{\frac{1}{p'_2}} \nonumber \\
&\lesssim F_R[u_1]^{\frac{1}{p_2}}\, R^{-2\sigma+ \frac{n+2\sigma}{p'_2}}\left(\int_{\R^n} \langle \tilde{x}\rangle^{-n-2\bar{\sigma}}\, d\tilde{x}\right)^{\frac{1}{p'_2}} \quad (\text{by }\eqref{key_Laplace})\nonumber \\
&\lesssim R^{-2\sigma+ \frac{n+2\sigma}{p'_2}}F_R[u_1]^{\frac{1}{p_2}}. \label{blow-up.t.4}
\end{align}
Thanks to the assumption (\ref{eta.Cond}), by the similar manner in terms of estimating $F_{2,R}[u_{\mathtt{k}}]$ one gets
\begin{equation}
\big|F_{4,R}[u_{\mathtt{k}}]\big|\lesssim R^{-4\sigma+ \frac{n+2\sigma}{p'_2}}F_R[u_1]^{\frac{1}{p_2}}. \label{blow-up.t.5}
\end{equation}
As a consequence, combining from (\ref{blow-up.t.1}) to (\ref{blow-up.t.5}) we may conclude
$$ F_R[u_{\mathtt{k}}] + \int_{\R^n} \big((-\Delta)^\sigma u_{01}(x)+ u_{01}(x)+ u_{11}(x)\big)\psi_R(x)\,dx \lesssim R^{-2\sigma+ \frac{n+2\sigma}{p'_2}}F_R[u_1]^{\frac{1}{p_2}}. $$
Analogously, we also arrive at the following estimates:
\begin{align}
F_R[u_1] &+ \int_{\R^n} \big((-\Delta)^\sigma u_{02}(x)+ u_{02}(x)+ u_{12}(x)\big)\psi_R(x)\,dx \lesssim R^{-2\sigma+ \frac{n+2\sigma}{p'_3}}F_R[u_2]^{\frac{1}{p_3}}, \nonumber \\
F_R[u_2] &+ \int_{\R^n} \big((-\Delta)^\sigma u_{03}(x)+ u_{03}(x)+ u_{13}(x)\big)\psi_R(x)\,dx \lesssim R^{-2\sigma+ \frac{n+2\sigma}{p'_4}}F_R[u_3]^{\frac{1}{p_4}}, \nonumber \\
\vdots \nonumber \\
F_R[u_{\mathtt{k}-1}] &+ \int_{\R^n} \big((-\Delta)^\sigma u_{0{\mathtt{k}}}(x)+ u_{0{\mathtt{k}}}(x)+ u_{1{\mathtt{k}}}(x)\big)\psi_R(x)\,dx \lesssim R^{-2\sigma+ \frac{n+2\sigma}{p'_1}}F_R[u_{\mathtt{k}}]^{\frac{1}{p_1}}. \label{blow-up.t.6}
\end{align}
Thanks to the hypothesis (\ref{blow-up.1}), there exists a sufficiently large constant $R_0$ so that it holds for any $R>R_0$:
\begin{equation} \label{blow-up.Initial}
\int_{\R^n} \big((-\Delta)^\sigma u_{0\ell}(x)+ u_{0\ell}(x)+ u_{1\ell}(x)\big)\psi_R(x)\,dx> 0
\end{equation}
with $\ell=1,2,\cdots,\mathtt{k}$. Thus, it follows that
\begin{align*}
F_R[u_{\mathtt{k}}] &\lesssim R^{-2\sigma+ \frac{n+2\sigma}{p'_2}}F_R[u_1]^{\frac{1}{p_2}}, \\
F_R[u_1] &\lesssim R^{-2\sigma+ \frac{n+2\sigma}{p'_3}}F_R[u_2]^{\frac{1}{p_3}}, \\
\vdots \\
F_R[u_{\mathtt{k}-2}] &\lesssim R^{-2\sigma+ \frac{n+2\sigma}{p'_{\mathtt{k}}}}F_R[u_{\mathtt{k}-1}]^{\frac{1}{p_{\mathtt{k}}}}.
\end{align*}
Now we plug the above chain of estimates into (\ref{blow-up.t.6}) successively to achieve
\begin{align*}
&F_R[u_{\mathtt{k}-1}] + \int_{\R^n} \big((-\Delta)^\sigma u_{0{\mathtt{k}}}(x)+ u_{0{\mathtt{k}}}(x)+ u_{1{\mathtt{k}}}(x)\big)\psi_R(x)\,dx \\
&\qquad \lesssim R^{-2\sigma+ \frac{n+2\sigma}{p'_1}}\Big(R^{-2\sigma+ \frac{n+2\sigma}{p'_2}}F_R[u_1]^{\frac{1}{p_2}}\Big)^{\frac{1}{p_1}} \\
&\qquad \quad= R^{-2\sigma\big(1+ \frac{1}{p_1}\big)+ (n+2\sigma)\big(1-\frac{1}{p_1p_2}\big)}F_R[u_1]^{\frac{1}{p_1p_2}} \\
&\qquad \lesssim R^{-2\sigma- \frac{2\sigma}{p_1}+ (n+2\sigma)\big(1-\frac{1}{p_1p_2}\big)}\Big(R^{-2\sigma+ \frac{n+2\sigma}{p'_3}}F_R[u_2]^{\frac{1}{p_3}}\Big)^{\frac{1}{p_1p_2}} \\
&\qquad \quad= R^{-2\sigma\big(1+ \frac{1}{p_1}+ \frac{1}{p_1p_2}\big)+ (n+2\sigma)\big(1-\frac{1}{p_1p_2p_3}\big)}F_R[u_2]^{\frac{1}{p_1p_2p_3}} \\
&\qquad \quad \vdots \\
&\qquad \lesssim R^{-2\sigma\big(1+ \frac{1}{p_1}+ \frac{1}{p_1p_2}+\cdots+\frac{1}{p_1p_2\cdots p_{\mathtt{k}-1}}\big)+ (n+2\sigma)\big(1-\frac{1}{p_1p_2\cdots p_{\mathtt{k}}}\big)}F_R[u_{\mathtt{k}-1}]^{\frac{1}{p_1p_2\cdots p_{\mathtt{k}}}},
\end{align*}
which is equivalent to
\begin{align}
&\int_{\R^n} \big((-\Delta)^\sigma u_{0{\mathtt{k}}}(x)+ u_{0{\mathtt{k}}}(x)+ u_{1{\mathtt{k}}}(x)\big)\psi_R(x)\,dx \nonumber \\
&\qquad \lesssim R^{-2\sigma\big(1+ \frac{1}{p_1}+ \frac{1}{p_1p_2}+\cdots+\frac{1}{p_1p_2\cdots p_{\mathtt{k}-1}}\big)+ (n+2\sigma)\big(1-\frac{1}{p_1p_2\cdots p_{\mathtt{k}}}\big)}F_R[u_{\mathtt{k}-1}]^{\frac{1}{p_1p_2\cdots p_{\mathtt{k}}}}- F_R[u_{\mathtt{k}-1}]. \label{blow-up.t.7}
\end{align}
Moreover, the application of the inequality
$$ A\,y^\rho- y \le A^{\frac{1}{1-\rho}} \quad \text{ for any } A>0,\, y \ge 0 \text{ and } 0< \rho< 1 $$
implies immediately
\begin{align}
&R^{-2\sigma\big(1+ \frac{1}{p_1}+ \frac{1}{p_1p_2}+\cdots+\frac{1}{p_1p_2\cdots p_{\mathtt{k}-1}}\big)+ (n+2\sigma)\big(1-\frac{1}{p_1p_2\cdots p_{\mathtt{k}}}\big)}F_R[u_{\mathtt{k}-1}]^{\frac{1}{p_1p_2\cdots p_{\mathtt{k}}}}- F_R[u_{\mathtt{k}-1}] \nonumber \\
&\hspace{5cm} \le R^{-\frac{2\sigma(1+p_{\mathtt{k}}+\cdots+ p_2p_3\cdots p_{\mathtt{k}})}{p_1p_2\cdots p_{\mathtt{k}}-1}+ n} \label{blow-up.t.8}
\end{align}
for all $R > R_0$. So, collecting the estimates (\ref{blow-up.t.7}) and (\ref{blow-up.t.8}) one obtains
\begin{equation} \label{blow-up.t.9}
\int_{\R^n} \big((-\Delta)^\sigma u_{0{\mathtt{k}}}(x)+ u_{0{\mathtt{k}}}(x)+ u_{1{\mathtt{k}}}(x)\big)\psi_R(x)\,dx \le C R^{-\frac{2\sigma(1+p_{\mathtt{k}}+\cdots+ p_2p_3\cdots p_{\mathtt{k}})}{p_1p_2\cdots p_{\mathtt{k}}-1}+ n}
\end{equation}
for any $R > R_0$. It is clear that the assumption (\ref{blow-up.2}) is equivalent to
$$ -\frac{2\sigma(1+p_{\mathtt{k}}+\cdots+ p_2p_3\cdots p_{\mathtt{k}})}{p_1p_2\cdots p_{\mathtt{k}}-1}+ n< 0. $$
Then, letting $R \to \ity$ in (\ref{blow-up.t.9}) we obtain a contradiction to \eqref{blow-up.Initial}. Summarizing, the proof is completed.

%===================================================
\section{Lifespan estimates} \label{Sec.Lifespan}
This section is to devote to the study of the lifespan estimates for solutions to \eqref{MainSytem}. To establish this purpose, let us consider the initial data $\big(\varepsilon u_{0\ell}(x),\varepsilon u_{1\ell}(x)\big)$ instead of $\big(u_{0\ell}(x),u_{1\ell}(x)\big)$ for \eqref{MainSytem}, where $\varepsilon$ is a small positive constant which presents the size of initial data.

\subsection{Proof of Theorem \ref{Lifespan.Thoerem-1}} Assume that $u= u(t,x)$ is a local solution to \eqref{MainSytem} in $[0,T_\e)\times \R^n$ with $T_\e= {\rm LifeSpan}(u)$. Then, repeating some arguments to get the estimate \eqref{blow-up.t.9} in the proof of Theorem \ref{Blow-up.Thoerem} we may conclude the following estimate:
$$ c\varepsilon \le C R^{-\frac{2\sigma(1+p_{\mathtt{k}}+\cdots+ p_2p_3\cdots p_{\mathtt{k}})}{p_1p_2\cdots p_{\mathtt{k}}-1}+ n}, $$
where a positive constant $c$ is appropriately chosen to fulfill
$$ \int_{\R^n} \big((-\Delta)^\sigma u_{0\ell}(x)+ u_{0\ell}(x)+ u_{1\ell}(x)\big)\psi_R(x)\,dx> c > 0 $$
for any $R > R_0$ and $\ell=1,2,\cdots,\mathtt{k}$ by the help of the condition (\ref{blow-up.1}). By letting $R \to T_\e^{\frac{1}{2\sigma}}$ in the last inequality, one arrives at
$$ c\varepsilon \le T_\e^{-\frac{1+p_{\mathtt{k}}+\cdots+ p_2p_3\cdots p_{\mathtt{k}}}{p_1p_2\cdots p_{\mathtt{k}}-1}+ \frac{n}{2\sigma}}, $$
that is,
$$ T_\e \le C\varepsilon^{-\frac{1}{\gamma_{\mathtt{k}}-\frac{n}{2\sigma}}}, $$
what we wanted to prove. Hence, the proof of Theorem \ref{Lifespan.Thoerem-1} is completed.

\subsection{Proof of Theorem \ref{Lifespan.Thoerem-2}} With $\ell= 1,2,\cdots,\mathtt{k}$ we introduce the following evolution spaces:
$$ X_\ell(t):= \mathcal{C}\big([0,t],H^\sigma\big) $$
endowed with the corresponding norms as follows:
$$ \|u_\ell\|_{X_\ell(t)}:= \sup_{0\le \tau \le t} \left((1+\tau)^{\frac{n}{4\sigma}-\alpha_\ell}\|u_\ell(\tau,\cdot)\|_{L^2}+ (1+\tau)^{\frac{n}{4\sigma}+\frac{1}{2}-\alpha_\ell}\big\||D|^{\sigma} u_\ell(\tau,\cdot)\big\|_{L^2}\right) $$
for $\ell= 1,2,\cdots,\mathtt{k}-1$, and
$$ \|u_\mathtt{k}\|_{X_\mathtt{k}(t)}:=  \sup_{0\le \tau \le t} \left((1+\tau)^{\frac{n}{4\sigma}}\|u_{\mathtt{k}}(\tau,\cdot)\|_{L^2}+ (1+\tau)^{\frac{n}{4\sigma}+\frac{1}{2}}\big\||D|^{\sigma} u_{\mathtt{k}}(\tau,\cdot)\big\|_{L^2}\right). $$
Here these parameters $\alpha_\ell$ with $\ell= 1,2,\cdots,\mathtt{k}-1$ are defined by
$$ \alpha_\ell= \begin{cases}
1-(p_1-1)\gamma_{\mathtt{k}}>0 &\text{ if } \ell=1, \\
1-(p_\ell-1)\gamma_{\mathtt{k}}+p_\ell\alpha_{\ell-1} &\text{ if } \ell=2,3,\cdots \mathtt{k}-1.
\end{cases} $$
Now let us introduce the solution space $X(t)$ to the system \eqref{MainSytem} in the form
$$ X(t)= X_1(t) \times X_2(t) \times \cdots \times X_\mathtt{k}(t) $$
together with the norm of $u= (u_1,u_2,\cdots,u_{\mathtt{k}})$ as
$$ \|u\|_{X(t)}= \|(u_1,u_2,\cdots,u_{\mathtt{k}})\|_{X(t)}:= \sum_{\ell=1}^{\mathtt{k}}\|u_\ell\|_{X_\ell(t)}. $$
To establish the lower bound estimates for the lifespan, it is necessary to indicate the series of the following inequalities:
\begin{align}
\|u_1\|_{X_1(t)} &\le C_0\varepsilon +C^{u_1}_1(1+t)^{1-\frac{n}{2\sigma}(p_1-1)-\alpha_1}\|u_\mathtt{k}\|_{X_\mathtt{k}(t)}^{p_1}, \label{equation3.4} \\
\|u_\ell\|_{X_\ell(t)} &\le C_0\varepsilon +C^{u_\ell}_1(1+t)^{-\frac{n}{2\sigma}(p_\ell-1)+(p_\ell-1)\gamma_\mathtt{k}}\|u_{\ell-1}\|_{X_{\ell-1}(t)}^{p_\ell}, \label{equation3.5}
\end{align}
for $\ell= 2,3\cdots,\mathtt{k}$. Here $C_0$ stands for a suitable constant only dependent on $n$ and the initial data, while $C_1^{u_\ell}$ with $\ell= 1,2,\cdots,\mathtt{k}$ are positive constants independent of $t$. \medskip

Indeed, recalling the presentation formula of solutions \eqref{Sol.Rep-1} we may immediately arrive at the following estimate from the definition of the norm of $X(t)$ and Proposition \ref{prop.1.4.1}:
$$ \|u^{\rm ln}\|_{X(t)}\le \varepsilon C_0\big(n,\mathcal{I}[(u_{01},u_{11}),(u_{02},u_{12}),\cdots,(u_{0\mathtt{k}},u_{1\mathtt{k}})]\big). $$
For this reason, we deduce the proof of the inequalities (\ref{equation3.4}) and (\ref{equation3.5}) by showing the alternatives as follows:
\begin{align}
\|u^{\rm nl}_1\|_{X_1(t)} &\le C^{u_1}_1(1+t)^{1-\frac{n}{2\sigma}(p_1-1)-\alpha_1}\|u_\mathtt{k}\|_{X_\mathtt{k}(t)}^{p_1}, \label{equation3.6} \\
\|u^{\rm nl}_\ell\|_{X_\ell(t)} &\le C^{u_\ell}_1 (1+t)^{-\frac{n}{2\sigma}(p_\ell-1)+(p_\ell-1)\gamma_\mathtt{k}}\|u_{\ell-1}\|_{X_{\ell-1}(t)}^{p_\ell}, \label{equation3.7}
\end{align}
for $\ell= 2,3,\cdots,\mathtt{k}$. The employment of the fractional Gagliardo-Nirenberg inequality from Proposition \ref{fractionalGagliardoNirenberg} leads to
\begin{align*}
\||u_\ell(\tau,\cdot)|^{p_{\ell+1}}\|_{L^1} &\lesssim(1+\tau)^{-\frac{n}{2\sigma}(p_{\ell+1}-1)+\alpha_\ell p_{\ell+1}}{\|u_\ell\|^{p_{\ell+1}}_{X_\ell(t)}}, \\
\||u_\ell(\tau,\cdot)|^{p_{\ell+1}}\|_{L^2} &\lesssim(1+\tau)^{-\frac{n}{2\sigma}(p_{\ell+1}-\frac{1}{2})+\alpha_\ell p_{\ell+1}}{\|u_\ell\|^{p_{\ell+1}}_{X_\ell(t)}}, \\
\||u_\mathtt{k}(\tau,\cdot)|^{p_1}\|_{L^1} &\lesssim(1+\tau)^{-\frac{n}{2\sigma}(p_1-1)}{\|u_\mathtt{k}\|^{p_1}_{X_\mathtt{k}(t)}}, \\
\||u_\mathtt{k}(\tau,\cdot)|^{p_1}\|_{L^2} &\lesssim(1+\tau)^{-\frac{n}{2\sigma}(p_1-\frac{1}{2})}{\|u_\mathtt{k}\|^{p_1}_{X_\mathtt{k}(t)}},
\end{align*}
for any $0\le \tau\le t$ and $\ell= 1,2,\cdots,\mathtt{k}-1$, provided that the conditions \eqref{GN1A1} and \eqref{GN1A1} are true. Since the condition (\ref{blow-up.2}), it follows that
$$ -\frac{n}{2\sigma}(p_\mathtt{k}-1)+\alpha_{\mathtt{k}-1} p_\mathtt{k} <-1. $$
Performing the same procedures to the proof of Theorem \ref{Global-existence.Thoerem} one gets
\begin{align*}
\|u^{\rm nl}_\ell(\tau,\cdot)\|_{L^2} &\lesssim (1+\tau)^{-\frac{n}{4\sigma}+1-\frac{n}{2\sigma}(p_\ell-1)+\alpha_{\ell-1}p_\ell}\|u_{\ell-1}\|_{X_{\ell-1}(t)}^{p_\ell}, \\
\|u^{\rm nl}_1(\tau,\cdot)\|_{L^2} &\lesssim (1+\tau)^{-\frac{n}{4\sigma}+1-\frac{n}{2\sigma}(p_1-1)}\|u_{\mathtt{k}}\|_{X_{\mathtt{k}}(t)}^{p_1},
\end{align*}
i.e.
\begin{align*}
(1+\tau)^{\frac{n}{4\sigma}-\alpha_\ell}\|u^{\rm nl}_\ell(\tau,\cdot)\|_{L^2} &\lesssim (1+\tau)^{-\frac{n}{2\sigma}(p_\ell-1)+(p_\ell-1)\gamma_\mathtt{k}}\|u_{\ell-1}\|_{X_{\ell-1}(t)}^{p_\ell}, \\
(1+\tau)^{\frac{n}{4\sigma}-\alpha_1}\|u^{\rm nl}_1(\tau,\cdot)\|_{L^2} &\lesssim (1+\tau)^{1-\frac{n}{2\sigma}(p_1-1)-\alpha_1}\|u_\mathtt{k}\|_{X_\mathtt{k}(t)}^{p_1},
\end{align*}
for $\ell= 2,3,\cdots,\mathtt{k}$. Similarly, we also demonstrate the estimates
$$ (1+\tau)^{\frac{n}{4\sigma}-\alpha_\ell}\||D|^\sigma u^{\rm nl}_\ell(\tau,\cdot)\|_{L^2}\ \lesssim (1+\tau)^{-\frac{n}{2\sigma}(p_\ell-1)+(p_\ell-1)\gamma_\mathtt{k}}\|u_{\ell-1}\|_{X_{\ell-1}(t)}^{p_\ell} $$
for $\ell= 2,3,\cdots,\mathtt{k}$ and
$$ (1+\tau)^{\frac{n}{4\sigma}-\alpha_1}\||D|^\sigma u^{\rm nl}_1(\tau,\cdot)\|_{L^2}\ \lesssim (1+\tau)^{1-\frac{n}{2\sigma}(p_1-1)-\alpha_1}\|u_\mathtt{k}\|_{X_\mathtt{k}(t)}^{p_1}. $$
By using these obtained estimates linked to the definition of the corresponding norms in $X_\ell(t)$ with $\ell=1,2,\cdots,\mathtt{k}$ we conclude immediately (\ref{equation3.4}) and (\ref{equation3.5}).

What's more, let us now assume that
$$ T= \sup \big\{t\in[0,T_\e)\,\, \text{ such that }\,\, G(t)= \|u\|_{X(t)}\le M\varepsilon \big\}, $$
where $M>0$ is a sufficiently large constant. In what follows, we denote by $C_{i,j}$ with $i,j\in\N$, suitable positive constants independent of $\varepsilon$. Because of $\|u_1\|_{X_1(T)}\le \|u\|_{X(T)}\le M\varepsilon$, it follows successively from (\ref{equation3.5}) that the following estimates hold:
\begin{align}
\|u_2\|_{X_2(T)} &\le C_0\varepsilon + C_1^{u_2}(1+T)^{\beta_2} M^{p_2}\varepsilon^{p_2}, \nonumber \\
\|u_3\|_{X_3(T)} &\le C_0\varepsilon + C_1^{u_3}(1+T)^{\beta_3}\big(C_0\varepsilon + C_1^{u_2}(1+T)^{\beta_2} M^{p_2}\varepsilon^{p_2}\big)^{p_3} \nonumber \\
&\le C_0\varepsilon + C_{2,3}(1+T)^{\beta_3}\varepsilon^{p_3}+C_{3,3}(1+T)^{\beta_2p_3+\beta_3}M^{p_2p_3}\varepsilon^{p_2p_3}, \nonumber \\
& \vdots \nonumber \\
\|u_\mathtt{k}\|_{X_\mathtt{k}(T)} &\le C_0\varepsilon + C_{2,\mathtt{k}}(1+T)^{\beta_\mathtt{k}}\varepsilon^{p_\mathtt{k}}+\cdots+C_{\mathtt{k},\mathtt{k}}(1+T)^{\beta_\mathtt{k}+\cdots+\beta_{\mathtt{k}-1}p_\mathtt{k}}M^{p_2p_3\cdots p_\mathtt{k}}\varepsilon^{p_2p_3\cdots p_\mathtt{k}}, \label{equation3.9}
\end{align}
where the notations $\beta_\ell$ with $\ell=2,3,\cdots,\mathtt{k}$ are given by
$$ \beta_\ell= -\frac{n}{2\sigma}(p_\ell-1)+(p_\ell-1)\gamma_\mathtt{k}. $$
Then, substituting the estimates (\ref{equation3.9}) into (\ref{equation3.6}) results
\begin{align*}
\|u_1\|_{X_1(T)} &\le C_0\varepsilon + C_{2,1}(1+T)^{1-\frac{n}{2\sigma}(p_1-1)-\alpha_1}\varepsilon ^{p_1}+C_{3,1}(1+T)^{1-\frac{n}{2\sigma}(p_1-1)-\alpha_1+\beta_\mathtt{k}p_1}\varepsilon ^{p_1p_\mathtt{k}} \\
&\qquad +\cdots+C_{\mathtt{k}+1,1}(1+T)^{1-\frac{n}{2\sigma}(p_1-1)-\alpha_1+\beta_\mathtt{k}p_1+\cdots+\beta_{\mathtt{k}-1}p_\mathtt{k}p_1}M^{p_1p_2p_3\cdots p_\mathtt{k}}\varepsilon^{p_1p_2p_3\cdots p_\mathtt{k}} \\
&\le \varepsilon \Big(C_0 + C_{2,1}(1+T)^{1-\frac{n}{2\sigma}(p_1-1)-\alpha_1}\varepsilon ^{p_1-1}+C_{3,1}(1+T)^{1-\frac{n}{2\sigma}(p_1-1)-\alpha_1+\beta_\mathtt{k}p_1}\varepsilon ^{p_1p_\mathtt{k}-1} \\
&\qquad +\cdots+C_{\mathtt{k}+1,1}(1+T)^{1-\frac{n}{2\sigma}(p_1-1)-\alpha_1+\beta_\mathtt{k}p_1+\cdots+\beta_{\mathtt{k}-1}p_\mathtt{k}p_1}M^{p_1p_2p_3\cdots p_\mathtt{k}}\varepsilon^{p_1p_2p_3\cdots p_\mathtt{k}-1}\Big).
\end{align*}
Let us take a sufficiently large constant $M>0$ such that $0<C_0<\frac{M}{4\mathtt{k}}$. Suppose that all the following estimates occur:
$$ 4\mathtt{k}C_{2,1}M^{-1}(1+T)^{1-\frac{n}{2\sigma}(p_1-1)-\alpha_1}\varepsilon ^{p_1-1}<1, $$
$$ 4\mathtt{k}C_{3,1}M^{-1}(1+T)^{1-\frac{n}{2\sigma}(p_1-1)-\alpha_1\beta_\mathtt{k}p_1}\varepsilon ^{p_1p_\mathtt{k}-1}<1, $$
$$ \vdots$$
$$ 4\mathtt{k}C_{\mathtt{k}+1,1}(1+T)^{1-\frac{n}{2\sigma}(p_1-1)-\alpha_1\beta_\mathtt{k}p_1+\cdots+\beta_{\mathtt{k}-1}p_\mathtt{k}p_1}M^{p_1p_2p_3\cdots p_\mathtt{k}-1}\varepsilon ^{p_1p_2p_3\cdots p_\mathtt{k}-1}<1. $$
Consequently, we obtain
\begin{equation}\label{LS.Proof-1}
\|u_1\|_{X_1(T)}<\f{\mathtt{k}+1}{4\mathtt{k}}M\varepsilon.
\end{equation}
Furthermore, if we assume
$$ 4\mathtt{k}C_1^{u_\ell}(1+T)^{\beta_\ell}M^{p_\ell-1}\varepsilon^{p_\ell-1}<1 $$
for any $\ell=2,3,\cdots,\mathtt{k}$, then it yields from \eqref{equation3.5} that
\begin{equation} \label{LS.Proof-2}
\|u_\ell\|_{X_\ell(T)} \le \varepsilon \big(C_0+C_1^{u_\ell}(1+T)^{\beta_\ell}M^{p_\ell}\varepsilon^{p_\ell}\big)< \frac{1}{2\mathtt{k}}M\varepsilon,
\end{equation}
where we have used the relation $\|u_{\ell-1}\|_{X_{\ell-1}(T)}\le \|u\|_{X(T)}\le M\varepsilon$. Combining both the estimates \eqref{LS.Proof-1} and \eqref{LS.Proof-2} we have
\begin{equation} \label{equation3.10}
G(T)=\|u\|_{X(T)}=\sum_{\ell=1}^{\mathtt{k}}\|u_\ell\|_{X_\ell(T)}<\frac{3\mathtt{k}-1}{4\mathtt{k}}M\varepsilon <M\varepsilon.
\end{equation}
Observing that $G(t)$ is a continuous function for any $t\in(0,T_\e)$. Hence, it implies from (\ref{equation3.10}) that there exists a value $T_0\in(T,T_\e)$ enjoying $G(T_0)\le M\varepsilon$. This is a contradiction to the definition of $T$ which we introduced before. Therefore, we claim that one of the following estimates are true:
$$ 4\mathtt{k}C_{2,1}M^{-1}(1+T)^{1-\frac{n}{2\sigma}(p_1-1)-\alpha_1}\varepsilon ^{p_1-1}\geq1, $$
$$ 4\mathtt{k}C_{3,1}M^{-1}(1+T)^{1-\frac{n}{2\sigma}(p_1-1)-\alpha_1\beta_\mathtt{k}p_1}\varepsilon ^{p_1p_\mathtt{k}-1}\geq1, $$
$$ \vdots $$
$$ 4\mathtt{k}C_{{\mathtt{k}+1},1}(1+T)^{1-\frac{n}{2\sigma}(p_1-1)-\alpha_1\beta_\mathtt{k}p_1+\cdots+\beta_{\mathtt{k}-1}p_\mathtt{k}p_1}M^{p_1p_2p_3\cdots p_\mathtt{k}-1}\varepsilon ^{p_1p_2p_3\cdots p_\mathtt{k}-1}\geq1, $$
$$ 4\mathtt{k}C_1^{u_\ell}(1+T)^{\beta_\ell}M^{p_\ell-1}\varepsilon^{p_\ell-1}\geq1, $$
for any $\ell=2,3,\cdots,\mathtt{k}$. In this way, one catches the blow-up time $T_\e$ satisfying
\begin{align*}
T_\e &\geq c\varepsilon^{-\min\left\{\frac{p_1-1}{\beta_1},\cdots,\frac{p_\mathtt{k}-1}{\beta_\mathtt{k}},\cdots,\frac{p_1p_2p_3\cdots p_\mathtt{k}-1}{1-\frac{n}{2\sigma}(p_1-1)-\alpha_1+\beta_\mathtt{k}p_1+\cdots+\beta_{\mathtt{k}-1}p_\mathtt{k}p_1}\right\}} = c\varepsilon^{-\frac{1}{\gamma_\mathtt{k}-\frac{n}{2\sigma}}},
\end{align*}
where we notice that
\begin{align*}
\beta_1 &=1-\frac{n}{2\sigma}(p_1-1)-\alpha_1=(p_1-1)\left(\gamma_\mathtt{k}-\frac{n}{2\sigma}\right)>0, \\
\beta_\ell &=-\frac{n}{2\sigma}(p_\ell-1)+(p_\ell-1)\gamma_\mathtt{k}=(p_\ell-1)\left(\gamma_\mathtt{k}-\frac{n}{2\sigma}\right)>0,
\end{align*}
for any $\ell=2,3,\cdots,\mathtt{k}$, and
$$ 1-\frac{n}{2\sigma}(p_1-1)-\alpha_1+\beta_\mathtt{k}p_1+\cdots+\beta_{\mathtt{k}-1}p_\mathtt{k}p_1=(p_1p_2p_3\cdots p_\mathtt{k}-1)\left(\gamma_\mathtt{k}-\frac{n}{2\sigma}\right)>0 $$
thanks to the assumption \eqref{blow-up.2}. All in all, we have finished proving Theorem \ref{Lifespan.Thoerem-2}.

%===================================================

\section*{Acknowledgments}
This paper was completed with Dr. Tuan Anh Dao's patient guidance and his great help. This research of Tuan Anh Dao was partly supported by Vietnam Ministry of Education and Training under grant number B2023-BKA-06.

%=================================================================================={References}

\end{document}